\documentclass[12pt,lenq]{amsart}
\usepackage{amsmath}
\usepackage{amscd}
\usepackage{amssymb}
\usepackage{amsbsy}
\usepackage{amsfonts}
\usepackage{latexsym}
\usepackage{graphics}
\usepackage{amsmath,amscd,latexsym}
\usepackage{multirow}
\usepackage{array}
\usepackage{paralist}
\usepackage{titletoc}

\theoremstyle{plain}
\newtheorem{theorem}{Theorem}[section]
\newtheorem{proposition}[theorem]{Proposition}

\newtheorem{remark}[theorem]{Remark}
\newtheorem{definition}[theorem]{Definition}

\newtheorem{main theorem}[theorem]{Main Theorem}

\newtheorem{question}[theorem]{Question}

\newlength\savewidth

\newcommand{\interior}{\operatorname{int}}

\newcommand{\ZZ}{\mathbb{Z}}
\newcommand{\QQ}{\mathbb{Q}}
\newcommand{\RR}{\mathbb{R}}
\newcommand{\CC}{\mathbb{C}}
\newcommand{\HH}{\mathbb{H}}
\newcommand{\QQQ}{\hat{\mathbb{Q}}}
\newcommand{\RRR}{\hat{\mathbb{R}}}

\newcommand{\Conway}{\mbox{\boldmath$S$}^{2}}
\newcommand{\Conways}
{(\mbox{\boldmath$S$}^{2},\mbox{\boldmath$P$})}
\newcommand{\PP}{\mbox{\boldmath$P$}}
\newcommand{\PConway}{\mbox{\boldmath$S$}}
\newcommand{\ptorus}{\mbox{\boldmath$T$}}
\newcommand{\OO}{\mbox{\boldmath$O$}}

\newcommand{\PSL}{\mbox{$\mathrm{PSL}$}}
\newcommand{\SL}{\mbox{$\mathrm{SL}$}}
\newcommand{\tr}{\mbox{$\mathrm{tr}$}}
\newcommand{\Isom}{\mbox{$\mathrm{Isom}$}}

\newcommand{\DD}{\mathcal{D}}
\newcommand{\Riley}{\mathcal{R}}
\newcommand{\Discrete}{\mathcal{D}}
\newcommand{\RGPC}[2]{\Gamma({#1};{#2})}

\newcommand{\RGP}[1]{\Gamma_{#1}}
\newcommand{\curve}{\mathcal{C}}
\newcommand{\plamination}{\mathcal{PL}}
\newcommand{\einv}{\mathcal{E}}

\newcommand{\Hecke}{\mbox{$G$}}
\newcommand{\orbo}{\mbox{\boldmath$O$}}
\newcommand{\orbs}{\mbox{\boldmath$S$}}

\newcommand{\cfr}{\mbox{\boldmath$a$}}

\newcommand{\svert}{\,|\,}

\newcommand{\llangle}{\langle\langle}
\newcommand{\rrangle}{\rangle\rangle}

\begin{document}

\title{Simple loops on 2-bridge spheres in Heckoid orbifolds
for 2-bridge links}

\author{Donghi Lee}
\address{Department of Mathematics\\
Pusan National University \\
San-30 Jangjeon-Dong, Geumjung-Gu, Pusan, 609-735, Republic of Korea}
\email{donghi@pusan.ac.kr}

\author{Makoto Sakuma}
\address{Department of Mathematics\\
Graduate School of Science\\
Hiroshima University\\
Higashi-Hiroshima, 739-8526, Japan}
\email{sakuma@math.sci.hiroshima-u.ac.jp}

\subjclass[2010]{Primary 57M25, 20F06 \\
\indent {The first author was supported by Basic Science Research Program
through the National Research Foundation of Korea(NRF) funded
by the Ministry of Education, Science and Technology(2012009996).
The second author was supported
by JSPS Grants-in-Aid 22340013.}}


\begin{abstract}
Following Riley's work,
for each $2$-bridge link $K(r)$ of slope $r\in\QQ$
and an integer or a half-integer $n$ greater than $1$,
we introduce the {\it Heckoid orbifold $\orbs(r;n)$}
and the {\it Heckoid group $\Hecke(r;n)=\pi_1(\orbs(r;n))$ of
index $n$ for $K(r)$}.
When $n$ is an integer,
$\orbs(r;n)$ is called an {\it even} Heckoid orbifold;
in this case, the underlying space is the exterior of $K(r)$,
and the singular set is the lower tunnel of $K(r)$ with index $n$.
The main purpose of this note is to announce answers to
the following questions for even Heckoid orbifolds.
(1) For an essential simple loop on a $4$-punctured sphere $\PConway$
in $\orbs(r;n)$ determined by the $2$-bridge sphere of $K(r)$,
when is it null-homotopic in $\orbs(r;n)$?
(2) For two distinct essential simple loops
on $\PConway$, when are they homotopic in $\orbs(r;n)$?
We also announce applications of these results to
character varieties, McShane's identity, and
epimorphisms from $2$-bridge link groups onto Heckoid groups.
\end{abstract}
\maketitle


\section{Introduction}
\begin{figure}
 \centering
\includegraphics{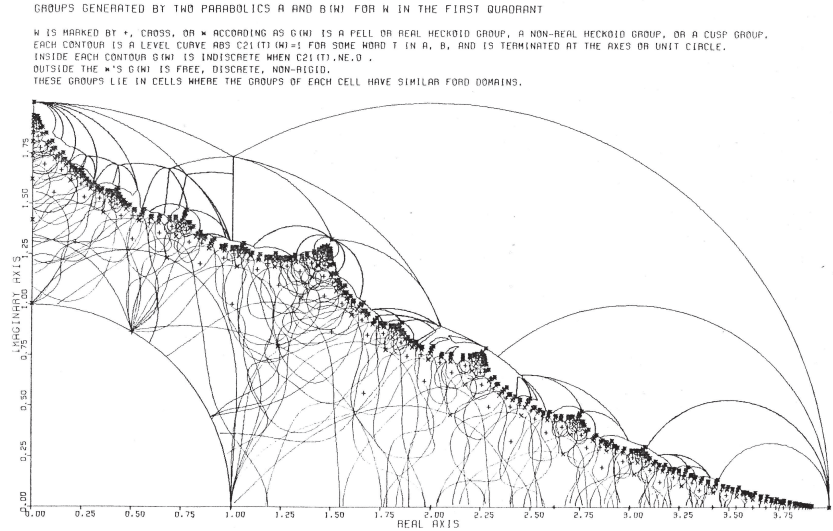}
\label{fig.Riley}
\caption{Riley's pioneering exploration of groups generated by two parabolic
transformations.
This specific copy of the picture
was obtained directly from Robert Riley by Masaaki Wada when he
visited SUNY Binghamton in February, 1991.
}
\end{figure}

In the late $70$'s, Robert Riley made a pioneering exploration of
groups generated by two parabolic transformations.
(His personal account of background history of
the exploration can be found in \cite{Riley3}.)
The computer-drawn Figure ~\ref{fig.Riley}
has been circulated among the experts
and has inspired many researchers in the fields of Kleinian groups and knot theory.
The diagram plots the set $\Discrete$ of those complex numbers $\omega$
in the first quadrant of the complex plane
for which the group
\[
G_{\omega}=
\left\langle
\begin{pmatrix}
	     1&1\\
	     0&1
	    \end{pmatrix},
\begin{pmatrix}
	     1&0\\
	     \omega&1
	    \end{pmatrix}
\right\rangle
\]
is discrete and non-free.

In \cite{Riley2}, in which he mentions
to the computer experiments
and describes the algebra behind the experiments,
Riley introduced an infinite collection of Laurent polynomials,
which he calls the Heckoid polynomials
associated with a $2$-bridge link $K$,
and observed
that these Heckoid polynomials define the affine representation varieties
of certain groups, the Heckoid groups for $K$.
He ``defines''
the Heckoid group of index $q\ge 3$ for $K$
to be the Kleinian group $G_{\omega}$,
where $\omega$ corresponds to
a ``right'' root of the Heckoid polynomials
(see \cite[the paragraph following Theorem ~A in p.390]{Riley2}).
The classical Hecke groups,
introduced in \cite{Hecke}, are essentially the simplest Heckoid groups.
He was convinced that
``the calculations for the diagram were pushed far enough to suggest
that all essential phenomena were depicted'',
and observed the following \cite[pp.391--392]{Riley2}:
\begin{enumerate}[\indent \rm (1)]
\item
Each $\omega\in \Discrete$ such that $G_{\omega}$ is torsion-free
corresponds to the hyperbolic structure on a $2$-bridge link complement.
\item
Each $\omega\in \Discrete$ such that $G_{\omega}$ is not torsion-free gives
a Heckoid group for some $2$-bridge link, or for the trivial knot.
\item
The points $\omega_q(K)$, $q=3,4,\dots$, determining the Heckoid groups
for a $2$-bridge link $K$ converges to a cusp $\omega_{\infty}(K)$
of the boundary of the region $\Riley$, the Riley slice of the Schottky space,
i.e., the set of complex numbers $\omega$
such that $G_{\omega}$ is discrete and free and
the quotient $\Omega(G_{\omega})/G_{\omega}$
of the domain of discontinuity $\Omega(G_{\omega})$ is homeomorphic to
the $4$-times punctured sphere.
\end{enumerate}

The classification of non-free $2$-parabolic generator Kleinian groups announced
by Agol ~\cite{Agol},
which generalizes the characterization of $2$-bridge links
obtained by \cite{Adams, Boileau-Zimmermann},
justifies Riley's anticipation.
In fact, the classification theorem says that
the Heckoid groups and the $2$-bridge link groups,
together with the $2$-parabolic generator subgroups of the classical Hecke groups
(which are Heckoid groups for the trivial knot),
are the only non-free Kleinian groups
generated by two parabolic transformations.

\begin{figure}
 \centering
 \includegraphics{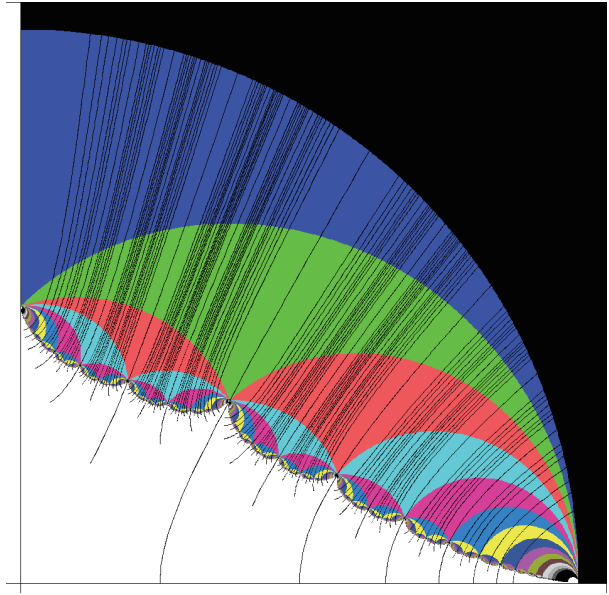}
 \caption{Riley slice of Schottky space together
  with pleating rays and their extensions}
 \label{fig.Yamashita}
\end{figure}

The announcement
made in the second author's joint work with
Akiyoshi, Wada and Yamashita in \cite[Section ~3 of Preface]{ASWY}
(cf. \cite[Remark ~6.1]{lee_sakuma_6})
identifies the points $\omega_q(K)$, $q=3,4,\dots$,
as discrete points in the extension of a rational pleating ray
studied by Keen and Series ~\cite{Keen-Series} (cf. \cite{Komori-Series}),
converging to rational boundary points of the Riley slice
(see Figure ~\ref{fig.Yamashita} copied from
\cite[Figure ~0.2b]{ASWY}).
Here, the endpoints of the extended rays correspond to $2$-bridge link groups.

This note is an announcement of the series of papers
~\cite{lee_sakuma_6, lee_sakuma_7, lee_sakuma_8, lee_sakuma_9}
on Heckoid groups motivated by Riley's work ~\cite{Riley2},
in particular his systematic construction
of epimorphisms from $2$-bridge link groups onto Heckoid groups
~\cite[the paragraph following Theorem ~B in p.391]{Riley2}.
The purposes of the series of papers are the following.
\begin{enumerate}[\indent \rm (1)]
\item
To give explicit combinatorial definitions of Heckoid groups
and Heckoid orbifolds for $2$-bridge links,
and to prove that the Heckoid groups are identified with
Kleinian groups generated by two parabolic transformations.
\item
To establish results for {\it even} Heckoid orbifolds
corresponding to those for $2$-bridge links
which are announced in \cite{lee_sakuma_0}
and proved in the preceding series of papers
\cite{lee_sakuma, lee_sakuma_2, lee_sakuma_3, lee_sakuma_4, lee_sakuma_5}.
\end{enumerate}
We note that the first purpose is essentially achieved in the announcement \cite{Agol},
but without proof.
To explain the second purpose, we
recall the preceding results mentioned in it.
They are concerned with
the following natural question,
which in turn is regarded as a variation of a question
raised by Minsky ~\cite[Question ~5.4]{Gordon}
on essential simple loops on Heegaard surfaces of $3$-manifolds.

\begin{question}
\label{question0}
{\rm
Let $K$ be a knot or a link in $S^3$
and $S$ a punctured sphere in the complement $S^3-K$
obtained from a bridge sphere of $K$.
\begin{enumerate}[\indent \rm (1)]
\item
Which essential simple loops on $S$ are null-homotopic or peripheral in
$S^3-K$?
\item
For two distinct essential simple loops on $S$,
when are they homotopic in $S^3-K$?
\end{enumerate}
}
\end{question}
In the preceding series of papers, we gave a complete answer to
the above question for $2$-bridge links
with applications to the study of character varieties, McShane's identity and
epimorphisms between $2$-bridge link groups.

In the new series of papers, we study the following variation of Question ~\ref{question0}
for $2$-bridge spheres of $2$-bridge links.

\begin{question} \label{question1}
{\rm
Consider the even Heckoid
orbifold $\orbs(r;n)$ of index $n$ for the $2$-bridge link $K(r)$,
and let $S$ be a $4$-holed sphere
in the underlying space $|\orbs(r;n)|$ obtained from
a $2$-bridge sphere of $K(r)$.
\begin{enumerate}[\indent \rm (1)]
\item
Which essential simple loops on $S$ are null-homotopic,
peripheral or torsion in $\orbs(r;n)$?
\item
For two distinct essential simple loops on $S$,
when are they homotopic in $\orbs(r;n)$?
\end{enumerate}
}
\end{question}

\begin{figure}[h]
\begin{center}
\includegraphics{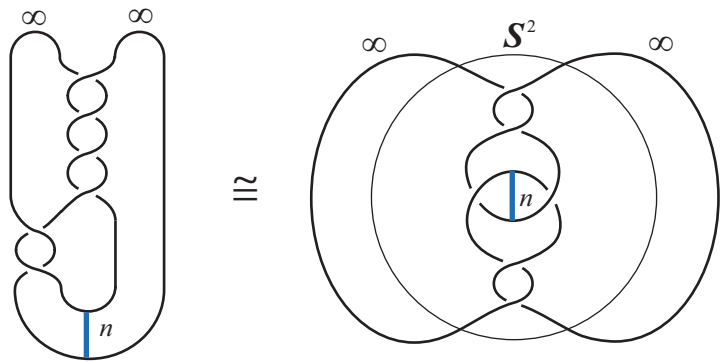}
\end{center}
\caption{
\label{fig.Hekoid-orbifold}
The Heckoid orbifold $\orbs(r;n)$.
The labels $\infty$
indicate the parabolic loci.
Here $(S^3,K(r))=(B^3,t(\infty))\cup (B^3,t(r))$ with $r=[4,2]=2/9$,
where $(B^3,t(r))$ and $(B^3,t(\infty))$, respectively,
are the inside and the outside of the bridge sphere $\Conway$.
The lower tunnel is
the core tunnel of $(B^3,t(r))$.}
\end{figure}

Here the
{\it even Heckoid orbifold, $\orbs(r;n)$,
of index $n$ for the $2$-bridge link $K(r)$} is
the $3$-orbifold as illustrated in Figure ~\ref{fig.Hekoid-orbifold}.
The underlying space
$|\orbs(r;n)|$ is the exterior, $E(K(r))=S^3-\interior N(K(r))$, of $K(r)$,
and the singular set is the lower tunnel of $K(r)$,
where the index of the singularity is $n$.

As applications, we give
(i) a systematic construction of upper-meridian-pair preserving
epimorphisms from $2$-bridge link groups onto Heckoid groups,
generalizing Riley's construction,
and (ii) the characterization of those epimorphisms onto even Heckoid groups.
We also announce applications
to character varieties and McShane's identity.

As in the preceding series of papers,
the key tool is small cancellation theory,
applied to two-generator and one-relator presentations
of even Heckoid groups.

This note is organized as follows.
In Section ~\ref{sec:Heckoid},
we explain combinatorial definitions of
Heckoid orbifolds and Heckoid groups
and explain that the Heckoid groups are identified with
$2$-parabolic generator Kleinian groups.
In Section ~\ref{statements},
we describe the main results.
In Sections ~\ref{application1} and \ref{application2},
we announce applications of the main results
to epimorphisms from $2$-bridge link groups onto Heckoid groups
and to character varieties and McShane's identity.
In Section ~\ref{group_presentation},
we give a brief explanation of the ideas for the proofs of the main results.

\section{Heckoid orbifold $\orbs(r;n)$ and Heckoid group $\Hecke(r;n)$}
\label{sec:Heckoid}

For a rational number $r \in \QQQ:=\QQ\cup\{\infty\}$,
let $K(r)$ be the $2$-bridge link of slope $r$,
which is defined as the sum
$(S^3,K(r))=(B^3,t(\infty))\cup (B^3,t(r))$
of rational tangles of slope $\infty$ and $r$
(cf. Figure ~\ref{fig.Hekoid-orbifold}).
The common boundary $\partial (B^3,t(\infty))= \partial (B^3,t(r))$
of the rational tangles is identified
with the {\it Conway sphere} $\Conways:=(\RR^2,\ZZ^2)/H$,
where $H$ is the group of isometries
of the Euclidean plane $\RR^2$
generated by the $\pi$-rotations around
the points in the lattice $\ZZ^2$.
Let $\PConway$ be the $4$-punctured sphere $\Conway-\PP$
in the link complement $S^3-K(r)$.
Any essential simple loop in $\PConway$, up to isotopy,
is obtained as
the image of a line of slope $s\in\QQQ$ in $\RR^2-\ZZ^2$
by the covering projection onto $\PConway$.
The (unoriented) essential simple loop in $\PConway$ so obtained
is denoted by $\alpha_s$.
We also denote by $\alpha_s$ the conjugacy class of
an element of $\pi_1(\PConway)$
represented by (a suitably oriented) $\alpha_s$.
The loops $\alpha_{\infty}$ and $\alpha_r$ bound disks
in $B^3-t(\infty)$ and $B^3-t(r)$, respectively.
Thus the {\it link group} $G(K(r))=\pi_1(S^3-K(r))$ is obtained as follows:
\[
G(K(r))=\pi_1(S^3-K(r))
\cong \pi_1(\PConway)/ \llangle\alpha_{\infty},\alpha_r\rrangle
\cong \pi_1(B^3-t(\infty))/\llangle\alpha_r\rrangle.
\]

For each rational number $r$ and an integer $n\ge 2$,
the {\it even Heckoid orbifold of index $n$ for the $2$-bridge link $K(r)$}
is the $3$-orbifold $\orbs(r;n)$,
such that the underlying space
$|\orbs(r;n)|$ is the exterior, $E(K(r))=S^3-\interior N(K(r))$, of $K(r)$,
and that the singular set is the lower tunnel of $K(r)$
(i.e., the core tunnel of $(B^3,t(\infty))$
in the sense of \cite[p.360]{lee_sakuma}),
where the index of the singularity is $n$
(see Figure ~\ref{fig.Hekoid-orbifold}).
We call the orbifold fundamental group $\pi_1(\orbs(r;n))$
the {\it Heckoid group of index $n$ for $K(r)$},
and denote it by $\Hecke(r;n)$.
Since the loop $\alpha_{r}$ is isotopic
to a meridional loop around the lower tunnel,
the even Hekoid group $\Hecke(r;n)=\pi_1(\orbs(r;n))$ ($n\ge 2$)
is obtained as follows:
\[
\Hecke(r;n)
=\pi_1(\orbs(r;n))
\cong\pi_1(\PConway)/ \llangle\alpha_{\infty},\alpha_r^n\rrangle
\cong \pi_1(B^3-t(\infty))/\llangle\alpha_r^n\rrangle.
\]

\begin{figure}[htbp]
\begin{center}
\includegraphics{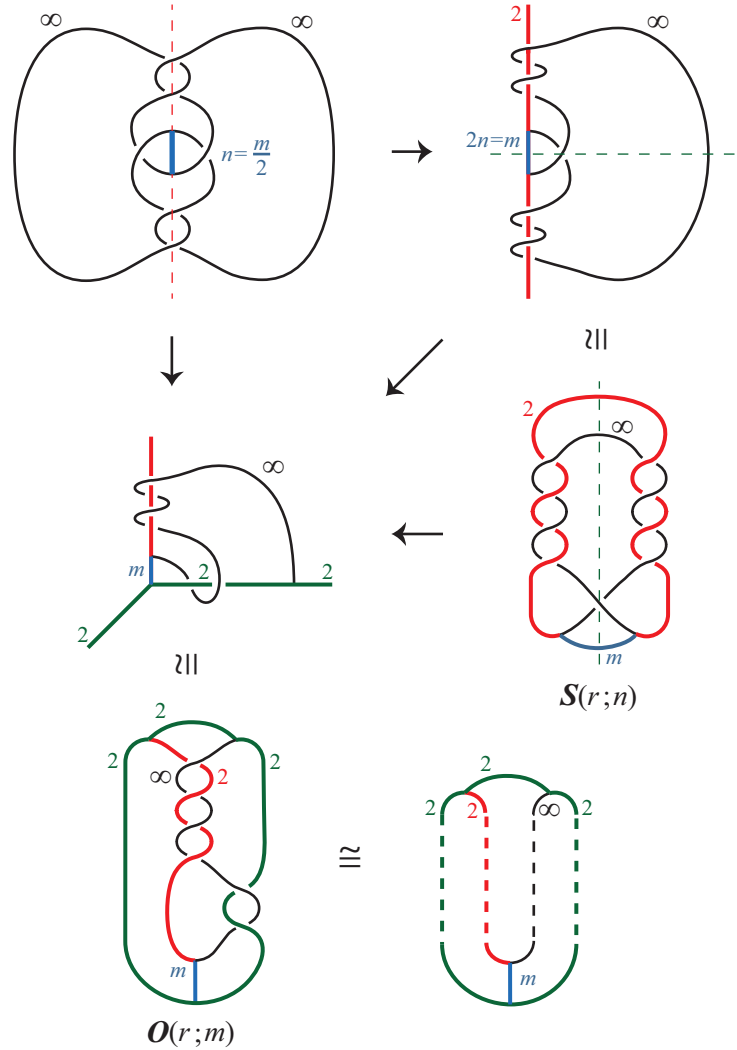}
\end{center}
\caption{\label{odd-Heckoid-orbifold1}
The case when $K(r)$ is a knot
and $m=2n>1$ is an odd integer.
Here $r=2/9=[4,2]$.
The odd Heckoid orbifold $\orbs(r;n)$ (middle right)
is a $\ZZ/2\ZZ$-covering of $\orbo(r;m)$ (lower left).
The upper left figure is not an orbifold, but
is a hyperbolic cone manifold.
The odd Heckoid orbifold $\orbs(r;n)$
is the quotient of the cone manifold
by the $\pi$-rotation around the axis containing the singular set.}
\end{figure}

The announcement by Agol ~\cite{Agol}
and the announcement
made in the second author's joint work with
Akiyoshi, Wada and Yamashita in \cite[Section ~3 of Preface]{ASWY}
suggest that the group $\Hecke(r;n)$ makes sense
even when $n$ is a half-integer greater than $1$.
We refer to \cite[Definition ~3.2]{lee_sakuma_6} for the definition of
the group $\Hecke(r;n)$ and the corresponding orbifold
$\orbs(r;n)$ when $n$ is a non-integral half-integer greater than $1$.
Roughly speaking,
$\orbs(r;n)$ is defined to be a $\ZZ/2\ZZ$-covering
of a certain orbifold $\orbo(r;m)$, with $m=2n$,
which is obtained from the quotient of $K(r)$
by the natural $(\ZZ/2\ZZ)^2$-symmetry
(see Figure ~\ref{odd-Heckoid-orbifold1} for the case when
$K(r)$ is a knot).
We call them the {\it odd Heckoid orbifold} and
the {\it odd Heckoid group}, respectively,
of index $n$ for $K(r)$.
A topological description of an odd Heckoid orbifold
is given by
\cite[Proposition ~5.3 and Figures ~5 and 6]{lee_sakuma_6}.
We note that similar orbifolds are studied by \cite{Jones-Reid, Mecchia-Zimmermann}.

\begin{remark}
{
\rm
Our terminology is slightly different from that of Riley ~\cite{Riley2},
where $\Hecke(r;n)$ is
called the Heckoid group of index ``$m$'' for $K(r)$ with $m=2n$.
The Heckoid orbifold $\orbs(r;n)$ and the Heckoid group $\Hecke(r;n)$
are {\it even} or {\it odd} according to whether Riley's index $m=2n$
is even or odd.
}
\end{remark}

The following theorem
was anticipated in \cite{Riley2}
and is contained in \cite{Agol} without proof.

\begin{theorem}
\label{thm.Kleinian_heckoid}
For a rational number $r$ and an integer or a half-integer $n>1$,
the Heckoid group $\Hecke(r;n)$ is isomorphic to a
geometrically finite Kleinian group
generated by two parabolic transformations.
\end{theorem}

A proof of this theorem is given in
\cite[Section ~6]{lee_sakuma_6}
by using the orbifold theorem
for pared orbifolds
\cite[Theorem ~8.3.9]{Boileau-Porti}
(cf. \cite{CHK, Kapovich}).
As noted in \cite{Agol},
the proof is analogous to the arguments in
\cite[Proof of Theorem ~9]{Jones-Reid}.

By this theorem and the topological description of
odd Heckoid orbifolds
(\cite[Proposition ~5.3]{lee_sakuma_6}),
we obtain the following proposition,
which shows a significant difference
between odd and even Heckoid groups
(see \cite[Section ~6]{lee_sakuma_6}).

\begin{proposition}
\label{prop-not-one-relator}
Any odd Heckoid group is not a one-relator group.
\end{proposition}

\section{Main results}
\label{statements}

Let $\DD$ be the {\it Farey tessellation}, that is,
the tessellation of the upper half
space $\HH^2$ by ideal triangles which are obtained
from the ideal triangle with the ideal vertices $0, 1, \infty \in \QQQ$
by repeated reflection in the edges.
Then $\QQQ$ is identified with the set of the ideal vertices of $\DD$.
For each $r\in \QQQ$, let $\RGP{r}$ be the group of automorphisms of
$\DD$ generated by reflections in the edges of $\DD$
with an endpoint $r$.
It should be noted that $\RGP{r}$
is isomorphic to the infinite dihedral group
and that the region bounded by two adjacent edges of $\DD$
with an endpoint $r$ is a fundamental domain
for the action of $\RGP{r}$ on $\HH^2$.
For an integer $m$,
let $C_r(m)$ be the group of automorphisms of $\DD$
generated by the parabolic transformation,
centered on the vertex $r$, by $m$ units
in the clockwise direction.

\begin{figure}[h]
\begin{center}
\includegraphics{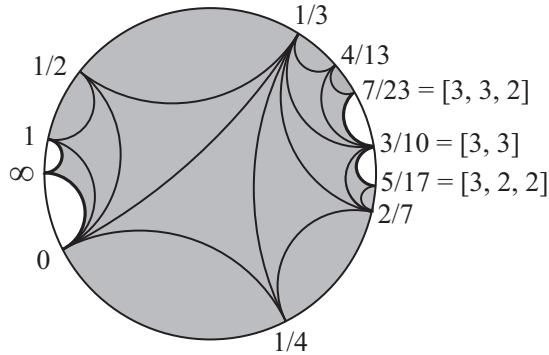}
\end{center}
\caption{\label{fig.fd_orbifold}
A fundamental domain of $\RGPC{r}{n}$ in the
Farey tessellation (the shaded domain) for $r=3/10=${\scriptsize $\cfrac{1}{3+\cfrac{1}{3}}$}\,$=:[3,3]$ and $n=2$.}
\end{figure}

For $r$ a rational number and $n$ an integer or a half-integer
greater than $1$,
let $\RGPC{r}{n}$ be the group generated by $\RGP{\infty}$ and $C_r(2n)$.
Suppose that $r$ is not an integer.
Then $\RGPC{r}{n}$ is the free product
$\RGP{\infty}*C_r(2n)$
having a fundamental domain, $R$, shown in Figure ~\ref{fig.fd_orbifold}.
Here, $R$ is obtained as the intersection of fundamental domains
for $\RGP{\infty}$ and $C_r(2n)$, and so
$R$ is bounded by the following two pairs of Farey edges:
\begin{enumerate}[\indent \rm (1)]
\item
the pair of adjacent Farey edges with an endpoint $\infty$
which cuts off a region in $\bar\HH^2$ containing $r$, and
\item
a pair of Farey edges with an endpoint $r$
which cuts off a region in $\bar\HH^2$ containing $\infty$, such that
one edge is the image of the other by a generator of $C_r(2n)$.
\end{enumerate}

Let $\bar I(n;r)$ be the union of two closed intervals in $\partial \HH^2=\RRR$
obtained as the intersection of the closure of $R$ with $\partial \HH^2$.
Note that there is a pair $\{r_1, r_2\}$
of boundary points of $\bar I(n;r)$
such that $r_2$ is the image of $r_1$ by a generator of $C_r(2n)$.
Set $I(n;r)=\bar I(n;r) -\{r_i\}$ with $i=1$ or $2$.
Note that
$I(n;r)$ is the disjoint union of a closed interval and a half-open interval,
except for the special case when $r\equiv \pm1/p \pmod{\ZZ}$.

The following theorem proved in \cite{lee_sakuma_7}
is the starting point of all the results which we announce in this note.

\begin{theorem}
\label{thm:fundametal_domain}
Suppose that $r$ is a non-integral rational number and
that $n$ is an integer or a half integer greater than $1$.
Then, for any $s\in\QQQ$, there is a unique rational number
$s_0\in I(r;n) \cup\{\infty, r\}$
such that $s$ is contained in the $\RGPC{r}{n}$-orbit of $s_0$.
Moreover $\alpha_s$ is homotopic to $\alpha_{s_0}$ in $\orbs(r;n)$.
In particular, if $s_0=\infty$, then $\alpha_s$ is
null-homotopic in $\orbs(r;n)$.
\end{theorem}

The following theorem proved in \cite{lee_sakuma_7}
gives an answer to
the first question in Question ~\ref{question1}(1).

\begin{theorem}
\label{thm:null-homotopy}
Suppose that $r$ is a non-integral rational number and
that $n$ is an integer with $n \ge 2$.
Then the loop $\alpha_s$ is null-homotopic in $\orbs(r;n)$
if and only if $s$ belongs to the $\RGPC{r}{n}$-orbit of $\infty$.
In other words, if $s\in I(r;n) \cup\{r\}$,
then $\alpha_s$ is not null-homotopic in $\orbs(r;n)$.
\end{theorem}

The following theorems, which we prove in
\cite{lee_sakuma_8},
give answers to the remainder of Question ~\ref{question1}(1)
and Question ~\ref{question1}(2).

\begin{theorem}
\label{thm:conjugacy}
Suppose that $r$ is a non-integral rational number and
that $n$ is an integer with $n \ge 2$.
For two rational numbers $s$ and $s'$,
the simple loops $\alpha_s$ and $\alpha_{s'}$ are homotopic in
$\orbs(r;n)$
if and only if $s$ and $s'$
belong to the same $\RGPC{r}{n}$-orbit.
In other words,
for distinct
$s, s'\in I(r;n) \cup\{\infty, r\}$,
the simple loops $\alpha_s$ and $\alpha_{s'}$ are not homotopic in
$\orbs(r;n)$.
\end{theorem}

\begin{theorem}
\label{thm:peripheral}
Suppose that $r$ is a non-integral rational number and
that $n$ is an integer with $n \ge 2$.
Then the following hold.
\begin{enumerate}[\indent \rm (1)]
\item
The loop $\alpha_s$ is peripheral in $\orbs(r;n)$
if and only if $s$ belongs to the $\RGPC{r}{n}$-orbit of $\infty$.
\item
The loop $\alpha_s$ is torsion in $\orbs(r;n)$
if and only if $s$ belongs to the $\RGPC{r}{n}$-orbit of $\infty$ or $r$.
\end{enumerate}
In other words, there is no rational number $s\in I(r;n)$
for which the simple loop $\alpha_s$ is peripheral or torsion in $\orbs(r;n)$.
\end{theorem}

In the above theorem,
we say that $\alpha_s$ is {\it peripheral} or {\it torsion}
if the conjugacy class $\alpha_s$ is represented by
a (possibly trivial) parabolic or elliptic transformation, respectively,
when we identify
$\Hecke(r;n)$ with
a Kleinian group generated by two parabolic transformations.

\section{Application to epimorphisms from $2$-bridge link groups onto Heckoid groups}
\label{application1}

In \cite{Riley2},
Riley discussed relations of the Heckoid polynomials with the polynomials
defining the nonabelian $\SL(2,\CC)$-representations
of $2$-bridge link groups introduced in \cite{Riley1},
and proved that each Heckoid polynomial divides
the nonabelian representation polynomials of $2$-bridge links $\tilde K$,
where $\tilde K$ belongs to an infinite collection
of $2$-bridge links determined by $K$ and index $q$.
This suggests that there are epimorphisms from the link group of $\tilde K$
onto the Heckoid group of index $q$ for $K$,
as observed in \cite[the paragraph following Theorem ~B in p.391]{Riley2}.
The following theorem proved in \cite{lee_sakuma_6}
gives a systematic construction of
epimorphisms from $2$-bridge link groups
onto Heckoid groups,
generalizing Riley's construction
(see \cite[Section 2]{lee_sakuma_6} for precise meaning).

\begin{theorem}
\label{thm:epimorophism}
Suppose that $r$ is a rational number and
that $n$ is an integer or a half-integer greater than $1$.
For $s\in\QQQ$, if $s$ or $s+1$ belongs to the $\RGPC{r}{n}$-orbit of $\infty$,
then there is an upper-meridian-pair-preserving epimorphism
from $G(K(s))$ to $\Hecke(r;n)$.
\end{theorem}

As noted in Remark ~\ref{rem:Riley's_theorem},
Riley's construction corresponds to the case
when $s$ belongs to the orbit of $\infty$
by the infinite cyclic subgroup $C_r(2n)$ of $\RGPC{r}{n}$.

For even Heckoid groups, the following converse to the above theorem
is proved in \cite{lee_sakuma_7}
as a corollary to Theorem ~\ref{thm:null-homotopy}
(see \cite[Section ~2]{lee_sakuma_7} for precise meaning).

\begin{theorem}
\label{thm:epimorophism2}
Suppose that $r$ is a non-integral rational number and
that $n$ is an integer with $n \ge 2$.
Then there is an upper-meridian-pair-preserving epimorphism
from $G(K(s))$ onto the even Heckoid group $\Hecke(r;n)$
if and only if $s$ or $s+1$ belongs to the
$\RGPC{r}{n}$-orbit of $\infty$.
\end{theorem}

We give a characterization of those rational numbers
which belong to the $\RGPC{r}{n}$-orbit of $\infty$.
Since $\Hecke(r;n)$ is isomorphic to $\Hecke(r+1;n)$,
we may assume in the remainder of this paper that $0<r<1$.
For the continued fraction expansion
\begin{center}
\begin{picture}(230,70)
\put(0,48){$\displaystyle{
r=[a_1,a_2, \dots,a_k]:=
\cfrac{1}{a_1+
\cfrac{1}{ \raisebox{-5pt}[0pt][0pt]{$a_2 \, + \, $}
\raisebox{-10pt}[0pt][0pt]{$\, \ddots \ $}
\raisebox{-12pt}[0pt][0pt]{$+ \, \cfrac{1}{a_k}$}
}}}$}
\end{picture}
\end{center}
where $k \ge 1$, $(a_1, \dots, a_k) \in (\mathbb{Z}_+)^k$, and $a_k \ge 2$,
let $\cfr$, $\cfr^{-1}$, $\epsilon\cfr$ and $\epsilon\cfr^{-1}$, with $\epsilon\in\{-,+\}$,
be the finite sequences defined as follows:
\begin{align*}
\cfr &= (a_1, a_2,\dots, a_k), \quad &
\cfr^{-1} &=(a_k,a_{k-1},\dots,a_1),\\
\epsilon\cfr &=(\epsilon a_1,\epsilon a_2,\dots,
\epsilon a_k), \quad &
\epsilon \cfr^{-1} &=(\epsilon a_k,\epsilon
a_{k-1},\dots,
\epsilon a_1).
\end{align*}
As stated in \cite{lee_sakuma_6},
the following proposition
can be proved by the argument in \cite[Section ~5.1]{Ohtsuki-Riley-Sakuma}.

\begin{proposition}
\label{prop:continued fraction}
Let $r$ be as above and $n$ an integer or a half-integer greater than $1$.
Set $m=2n$.
Then a rational number $s$
belongs to the $\RGPC{r}{n}$-orbit of $\infty$
if and only if
$s$ has the following continued
fraction expansion{\rm :}
\[
s=
2c+[\epsilon_1\cfr,mc_1,-\epsilon_1\cfr^{-1},
2c_2,\epsilon_2\cfr,mc_3,-\epsilon_2\cfr^{-1},
\dots,
2c_{2t-2},\epsilon_t \cfr,mc_{2t-1},-\epsilon_t \cfr^{(-1)}]
\]
for some positive integer $t$, $c\in\ZZ$,
$(\epsilon_1,\epsilon_2,\dots,\epsilon_t) \in \{-,+\}^t$
and $(c_1,c_2,\dots,c_{2t-1})\in\ZZ^{2t-1}$.
\end{proposition}

\begin{remark}
\label{rem:Riley's_theorem}
{\rm
Riley's Theorem ~B and Theorem ~3 in \cite{Riley2}
imply the following.
Let $\alpha$ and $\beta$ be relatively prime integers with $1\le \beta <\alpha$.
For integers $d\ge 2$, $m\ge 3$, and $e\ge 1$,
consider the $2$-bridge link $K(\beta^*/\alpha^*)$,
where $(\alpha^*,\beta^*)=(\alpha^d m,\alpha^{d-1}m(\alpha-\beta)+e)$.
Then there is an epimorphism from the link group $G(K(\beta^*/\alpha^*))$ onto
the Heckoid group $\Hecke(\beta/\alpha;n)$, where $n=m/2$.
This result corresponds to the case when $r=(\alpha-\beta)/\alpha$
and $s=[\cfr,mc,-\cfr^{-1}]$, where $c=\epsilon\alpha^{d}$ with $\epsilon=\pm 1$
in Proposition ~\ref{prop:continued fraction}.
In fact, a simple calculation shows
\[
s=(\alpha^{d-1}m(\alpha-\beta)+(-1)^k\epsilon)/(\alpha^d m)=\beta^*/\alpha^*,
\]
where $k$ is the length of $\cfr$ and $\epsilon$ is chosen so that
$(-1)^k\epsilon=e$.
Thus
Theorem ~\ref{thm:epimorophism} and
Proposition ~\ref{prop:continued fraction} imply that
there is an epimorphism from
the link group $G(K(\beta^*/\alpha^*))= G(K(s))$
onto the Heckoid group $\Hecke(r;n)\cong \Hecke(1-r;n)= \Hecke(\beta/\alpha;n)$,
recovering Riley's result.

The simplest case in Riley's construction is the case
where $(\alpha,\beta)=(3,1)$, $d=2$, $m=3$,
and so $(\alpha^*,\beta^*)=(27, 18\pm 1)$.
See \cite[Section ~1.3]{Riley3} for the important role
which was played by the $2$-bridge knot
$8_{11} \cong K(17/27)\cong K(19/27)$ in his work.
}
\end{remark}

\section{Application to character varieties and McShane's identity}
\label{application2}

For a rational number $r$ and an integer $n \ge 2$,
let $\rho_{r,n}$ be the $\PSL(2,\CC)$-representation of $\pi_1(\PConway)$
obtained as the composition
\[
\pi_1(\PConway) \to
\pi_1(\PConway)/ \llangle\alpha_{\infty},\alpha_r^n\rrangle
\cong
\Hecke(r;n)
\to
\Isom^+(\HH^3)
\cong
\PSL(2,\CC),
\]
where the last homomorphism is the holonomy representation
of the pared hyperbolic orbifold $\orbs(r;n)$.

Now, let $\ptorus$ be the once-punctured torus
obtained as the quotient
$(\RR^2-\ZZ^2)/\ZZ^2$,
and let $\OO$ be the orbifold
$(\RR^2-\ZZ^2)/\hat H$
where $\hat H$ is the group generated by $\pi$-rotations
around the points in $(\frac{1}{2}\ZZ)^2$.
Note that $\OO$ is the orbifold with
underlying space a once-punctured sphere
and with three cone points of cone angle $\pi$.
The surfaces $\ptorus$ and $\PConway$, respectively,
are $\ZZ/2\ZZ$-covering and $(\ZZ/2\ZZ)^2$-covering of $\OO$,
and hence their fundamental groups are identified
with subgroups
of the orbifold fundamental group $\pi_1(\OO)$
of indices $2$ and $4$, respectively.
The $\PSL(2,\CC)$-representation $\rho_{r,n}$ of $\pi_1(\PConway)$
extends, in a unique way, to that of $\pi_1(\OO)$
(see \cite[Proposition ~2.2]{ASWY}),
and so we obtain, in a unique way, a
$\PSL(2,\CC)$-representation of $\pi_1(\ptorus)$
by restriction.
We continue to denote it by $\rho_{r,n}$.
Note that $\rho_{r,n}:\pi_1(\ptorus) \to \PSL(2,\CC)$
is {\it type-preserving}, i.e.,
it satisfies the following conditions.
\begin{enumerate}[\indent \rm (1)]
\item
$\rho_{r,n}$ is irreducible, i.e.,
its image does not have a common fixed point on
$\partial \HH^3$.
\item
$\rho_{r,n}$ maps a peripheral element of $\pi_1(\ptorus)$
to a parabolic transformation.
\end{enumerate}

By extending the concept of a geometrically infinite end
of a Kleinian group,
Bowditch ~\cite{Bowditch2}
introduced the notion of end invariants of a type-preserving
$\PSL(2,\CC)$-representation of $\pi_1(\ptorus)$.
Tan, Wong and Zhang ~\cite{Tan_Wong_Zhang_6} (cf. \cite{Tan_Wong_Zhang_1})
extended this notion (with slight modification)
to an arbitrary $\PSL(2,\CC)$-representation of $\pi_1(\ptorus)$.
(To be precise, ~\cite{Bowditch2} and \cite{Tan_Wong_Zhang_6}
treat $\SL(2,\CC)$-representations.
However, the arguments work for $\PSL(2,\CC)$-representations.)
To recall the definition of end invariants,
let $\curve$ be the set of free homotopy classes
of essential simple loops on $\ptorus$.
Then
$\curve$ is identified with $\QQQ$, the vertex set of
the Farey tessellation $\DD$ by the rule
$s\mapsto \beta_s$.
The projective lamination space $\plamination$ of $\ptorus$
is then identified with $\RRR:=\RR\cup\{\infty\}$
and contains $\curve$ as the dense subset
of rational points.

\begin{definition}
\label{def_end_invariant}
{\rm
Let $\rho$ be an $\PSL(2,\CC)$-representation of $\pi_1(\ptorus)$.

(1) An element $X\in \plamination$ is an {\it end invariant}
of $\rho$ if there exists a sequence of distinct elements
$X_n\in\curve$ such that
\begin{enumerate}[\indent \rm (i)]
\item $X_n\to X$, and

\item $\{|\tr \rho(X_n)|\}_n$ is bounded from above.
\end{enumerate}

(2) $\einv(\rho)$ denotes the set of end invariants of $\rho$.
}
\end{definition}
In the above definition, it should be noted that
$|\tr \rho(X_n)|$ is well-defined
though $\tr \rho(X_n)$ is defined only up to sign.

Tan, Wong and Zhang ~\cite{Tan_Wong_Zhang_1, Tan_Wong_Zhang_6}
showed that $\einv(\rho)$
is a closed subset of $\plamination$ and
proved various interesting properties of $\einv(\rho)$,
including a characterization of
those representations $\rho$
with $\einv(\rho)=\emptyset$ or $\plamination$,
generalizing a result of Bowditch ~\cite{Bowditch2}.
They also proved that $\einv(\rho)$ is a Cantor set in various cases
and proposed interesting conjectures.
In \cite{lee_sakuma_5},
the authors proved that
if $\rho_r:\pi_1(T)\to \PSL(2,\CC)$ corresponds to
the faithful discrete representation of
a hyperbolic $2$-bridge link group $G(K(r))$,
then $\einv(\rho_r)$ is equal to the limit set
of a certain discrete subgroup of $\Isom(\HH^2)$
and hence is a Cantor set.
By using Theorems ~\ref{thm:conjugacy} and \ref{thm:peripheral},
we obtain the following theorem
which determines the set $\einv(\rho_{r,n})$.

\begin{theorem}
For a non-integral rational number $r$ and an integer $n \ge 2$,
the set $\einv(\rho_{r,n})$ of end invariants of $\rho_{r,n}$
is equal to the limit set $\Lambda(\RGPC{r}{n})$
of the group $\RGPC{r}{n}$.
\end{theorem}

In \cite{lee_sakuma_5},
the authors used the results obtained
by the series of papers
\cite{lee_sakuma, lee_sakuma_2, lee_sakuma_3, lee_sakuma_4}
to give a variation of McShane's identity,
which expresses the cusp shape of a hyperbolic $2$-bridge link
in terms of the complex translation lengths
of simple loops on the bridge sphere.
(The original McShane's identity
was given by \cite{McShane0, McShane}
and variations and applications were given by
\cite{AMS, AMS2, Bowditch1, Bowditch2, Mirzakhani,
Tan_Wong_Zhang_1, Tan_Wong_Zhang_2, Tan_Wong_Zhang_4,
Tan_Wong_Zhang_5, Tan_Wong_Zhang_6, Tan_Wong_Zhang_7}.)
Theorems ~\ref{thm:null-homotopy}, \ref{thm:conjugacy} and \ref{thm:peripheral},
together with a result of Bowditch ~\cite{Bowditch2},
enable us to give a variation
of McShane's identity for Heckoid orbifolds for $2$-bridge links.

\begin{theorem}
\label{thm:McShane}
For a non-integral rational number $r$ and an integer $n\ge 2$,
the following hold:
\[
2\sum_{s\in \mathrm{int}I(r;n)}
\frac{1}{1+e^{l_{\rho_{r,n}}(\beta_s)}}
+
\sum_{s\in\partial \bar I(r;n)}\frac{1}{1+e^{l_{\rho_{r,n}}(\beta_s)}}
=-1.
\]
\end{theorem}
In the above theorem,
$\beta_s$ denotes the simple loop on $\ptorus$ of slope $s$.
We can refine the above theorem to show that
certain partial sums in the above series
have geometric meanings.
For details, please see \cite{lee_sakuma_9}.

\section{Idea of the proofs of Theorems ~\ref{thm:null-homotopy},
\ref{thm:conjugacy} and \ref{thm:peripheral}
}
\label{group_presentation}

Let $\{a,b\}$ be the standard meridian generator pair of
$\pi_1(B^3-t(\infty))$
as described in
\cite[Section ~3]{lee_sakuma}
(see also \cite[Section ~5]{lee_sakuma_0}).
Then $\pi_1(B^3-t(\infty))$ is identified with the free group $F(a,b)$.
For the rational number $r=q/p$, where $p$ and $q$ are relatively
prime positive integers,
let $u_r$ be the word in $\{a,b\}$ obtained as follows.
(For a geometric description, see \cite[Section ~5]{lee_sakuma_0}.)
Set $\epsilon_i = (-1)^{\lfloor iq/p \rfloor}$,
where $\lfloor x \rfloor$ is the greatest integer not exceeding $x$.
\begin{enumerate}[\indent \rm (1)]
\item If $p$ is odd, then
\[u_{q/p}=a\hat{u}_{q/p}b^{(-1)^q}\hat{u}_{q/p}^{-1},\]
where
$\hat{u}_{q/p} = b^{\epsilon_1} a^{\epsilon_2} \cdots b^{\epsilon_{p-2}} a^{\epsilon_{p-1}}$.

\item If $p$ is even, then
\[u_{q/p}=a\hat{u}_{q/p}a^{-1}\hat{u}_{q/p}^{-1},\]
where
$\hat{u}_{q/p} = b^{\epsilon_1} a^{\epsilon_2} \cdots a^{\epsilon_{p-2}} b^{\epsilon_{p-1}}$.
\end{enumerate}
Then $u_r\in F(a,b)\cong\pi_1(B^3-t(\infty))$
is represented by the simple loop $\alpha_r$,
and we obtain the following two-generator one-relator presentation
of the Heckoid group:
\[
\Hecke(r;n)
\cong
\pi_1(\PConway)/ \llangle\alpha_{\infty},\alpha_r^n\rrangle
\cong\pi_1(B^3-t(\infty))/\llangle \alpha_r^n\rrangle
\cong \langle a, b \, | \, u_r^n \rangle.
\]
Theorems ~\ref{thm:null-homotopy}, \ref{thm:conjugacy} and \ref{thm:peripheral}
are proved by applying small cancellation theory to the above presentation.
In the remainder of this section,
we present key ideas in the proofs of these theorems,
by assuming notation in \cite[Sections ~5--8]{lee_sakuma_0}.

Suppose that $r$ is a rational number with $0 < r< 1$
and that $n$ is an integer with $n \ge 2$.
Let $R$ be the symmetrized subset of $F(a, b)$ generated
by the single relator $u_{r}^n$ of the presentation
$\Hecke(r;n)=\langle a, b \svert u_r^n \rangle$.
Then we have the following proposition,
which enables us to
apply small cancellation theory to our problem.

\begin{proposition}[{\cite[Proposition ~4.4 and Corollary ~4.9]{lee_sakuma_7}}]
\label{prop:small_cancellation_condition_heckoid}
The symmetrized set $R$ satisfies $C(4n)$ and $T(4)$.
Thus every reduced $R$-diagram is a $[4, 4n]$-map.
\end{proposition}

By using basic formulas in small cancellation theory
(see \cite[Theorem ~V.3.1]{lyndon_schupp}),
we obtain the following proposition.

\begin{proposition}
[{\cite[Proposition ~4.11]{lee_sakuma_7}}]
\label{prop:key}
Let $M$ be an arbitrary connected and simply-connected $[4,4n]$-map
such that there is no vertex of degree $3$ in $\partial M$.
Put
\[
\begin{aligned}
A= &\ \text{\rm the number of vertices $v$ in $\partial M$ such that $d_M(v)=2$;} \\
B= &\ \text{\rm the number of vertices $v$ in $\partial M$ such that $d_M(v) \ge 4$.}
\end{aligned}
\]
Then $A \ge (4n-3)B+4n$.
In particular, there are at least $4n-2$ consecutive vertices of degree $2$
on $\partial M$.
\end{proposition}

This proposition implies that
if $\alpha_s$ represents the trivial element of $\Hecke(r;n)$,
then the cyclic word $(u_s)$ contains a subword of the cyclic word $(u_r^{\pm n})$
which is a product of $4n-1$ pieces and is not a product of $4n-2$ pieces.
Theorem ~\ref{thm:null-homotopy} is proved
by an inductive argument, based on this fact,
on the length of the (positive) continued fraction expansion of $r$.

\begin{remark}
\rm
{In \cite[Theorem ~3]{Newman} (cf. \cite[Theorem ~IV.5.5]{lyndon_schupp}),
Newman gives a powerful theorem
for the word problem for one relator groups with torsion,
which implies that if a cyclically reduced word $w$
represents the trivial element in
$\Hecke(r;n)\cong \langle a,b \svert u_r^n\rangle$,
then the cyclic word $(w)$
contains a subword of the cyclic word $(u_r^{\pm n})$
of length greater than $(n-1)/n=1-1/n$ times the length of $u_r^n$.
On the other hand,
the above proof of Theorem ~\ref{thm:null-homotopy}
is based on the fact that if $u_s$
represents the trivial element in $\Hecke(r;n)$,
then the cyclic word $(u_s)$ contains a subword of
the cyclic word $(u_r^{\pm n})$
which is a product of $4n-1$ pieces
but is not a product of less than $4n-1$ pieces.
Since $u_r^{\pm n}$ is a product of $4n$-pieces,
and since every piece has length less than a half of the length of $u_r$
(see \cite[Lemma ~4.2]{lee_sakuma_7}),
this implies that the cyclic word $(u_s)$ contains a subword of
the cyclic word $(u_r^{\pm n})$ of length greater than $1-1/(2n)$
times the length of $u_r^n$.
This suggests that
the above result of Newman
would not be strong enough to establish
Theorem ~\ref{thm:null-homotopy}.
}
\end{remark}

In order to prove Theorems ~\ref{thm:conjugacy} and \ref{thm:peripheral},
we establish a structure theorem
which describes the possible shapes of annular diagrams
over the presentation
$\Hecke(r;n) \cong \langle a, b \, | \, u_r^n \rangle$,
as in the structure theorem
established in \cite[Theorem ~4.9 and Corollary ~4.11]{lee_sakuma_2}
(cf. \cite[Theorem ~9.1]{lee_sakuma_0})
for annular diagrams over the presentation
$G(K(r)) \cong \langle a, b \svert u_r \rangle$.

\begin{figure}[h]
\begin{center}
\includegraphics{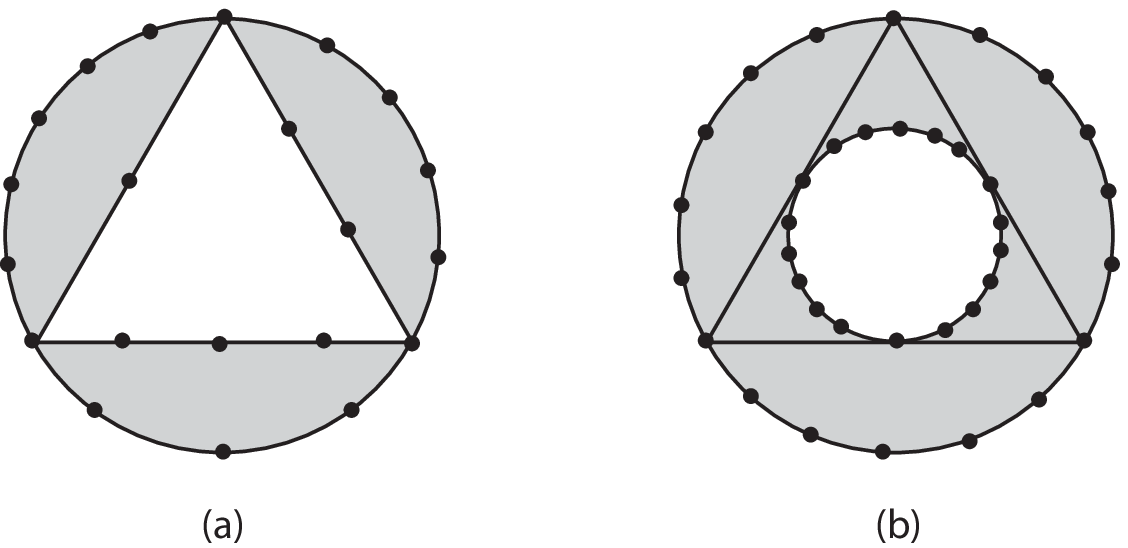}
\end{center}
\caption{\label{fig.layer}
$n=2$}
\end{figure}

\begin{theorem}[Structure Theorem]
\label{thm:annular_structure}
Let $M$ be an arbitrary nontrivial connected annular $[4,4n]$-map.
Suppose that there is no vertex of degree $1$ nor $3$ in $\partial M$
and that there are no $4n-2$ consecutive degree $2$ vertices in $\partial M$.
Let the outer and inner boundaries of $M$ be denoted by
$\sigma$ and $\tau$, respectively.
Then the following hold.
\begin{enumerate}[\indent \rm (1)]
\item The outer and inner boundaries $\sigma$ and $\tau$ are simple,
i.e., they are homeomorphic to the circle.

\item There is no edge contained in $\sigma \cap \tau$.

\item Every vertex of $M$ lies in $\partial M$.

\item $d_M(v)=2$ or $4$ for every vertex $v \in \partial M$.
Moreover, if $\sigma \cap \tau = \emptyset$,
then between any two vertices of degree $4$
there should occur exactly $(4n-3)$ vertices of degree $2$
on both $\sigma$ and $\tau$.

\item $d_M(D)=4n$ for every face $D \in M$.
\end{enumerate}
In particular,
Figure ~\ref{fig.layer} illustrates the only two possible shapes of $M$,
where the number of faces per layer is variable.
\end{theorem}

Theorems ~\ref{thm:conjugacy} and \ref{thm:peripheral} are proved
by an inductive argument, using this structure theorem,
on the length of the (positive) continued fraction expansion of $r$.
For details, please see \cite{lee_sakuma_8}.

\section*{Acknowledgements}
The second author would like to thank
(i) Gerhard Burde for drawing
his attention to the work of Riley ~\cite{Riley2} when he was staying in Frankfurt in $1997$,
(ii) Ian Agol for sending him a copy of the slides ~\cite{Agol} in $2002$, and
(iii) Gareth Jones for sending him Riley's personal account \cite{Riley3}
with a postscript by Brin, Jones and Singerman.
Both authors would like to thank Yeonhee Jang,
Makoto Ozawa and Toshio Saito
for helpful discussions and information
concerning the proof of Theorem ~\ref{thm.Kleinian_heckoid}.

\bibstyle{plain}

\end{document}